\documentclass[12pt]{elsarticle}
\usepackage[margin=1in]{geometry}

\usepackage{pgfplots}
\usepackage{tikz}
\usepackage{graphicx}
\graphicspath{ {./images/} }
\usepackage{amsmath,comment,enumitem}
\usepackage{mathrsfs}
\usepackage{amsthm,amssymb}
\usepackage{amsfonts}
\usepackage{mathtools}
\usepackage{hyperref}
\graphicspath{ {./images/} }
\usetikzlibrary{shapes}
\usepgfplotslibrary{polar}
\usetikzlibrary{decorations.markings}
\usetikzlibrary{backgrounds}
\pgfplotsset{every axis/.append style={
                    axis x line=middle,    % put the x axis in the middle
                    axis y line=middle,    % put the y axis in the middle
                    axis line style={<->,color=blue}, % arrows on the axis
                    xlabel={$x$},          % default put x on x-axis
                    ylabel={$y$},          % default put y on y-axis
            }}
\usepackage[utf8]{inputenc}

\newtheorem{theorem}{Theorem}
\newtheorem{lemma}[theorem]{Lemma}
\newtheorem{proposition}[theorem]{Proposition}
\newtheorem{corollary}[theorem]{Corollary}

\newtheorem{definition}{Definition}

\theoremstyle{remark}

\newtheorem{example}{Example}

\newtheorem{remark}[example]{Remark}

\begin{document}
\begin{frontmatter}

\author[add1,add2]{Mike Miller Eismeier\corref{cor1}}
\ead{Mike.Miller-Eismeier@uvm.edu}
\author[add1]{Aiden Sagerman}
\ead{afs2179@columbia.edu}
\address[add1]{Department of Mathematics, Columbia University}
\address[add2]{Department of Mathemstics \& Statistics, University of Vermont}
\cortext[cor1]{Corresponding author}

\title{Hyperplanes in abelian groups and twisted signatures}

\begin{abstract}
We investigate the following question: if $A$ and $A'$ are products of finite cyclic groups, when does there exist an isomorphism $f: A \to A'$ which preserves the union of coordinate hyperplanes (equivalently, so that $f(x)$ has some coordinate zero if and only if $x$ has some coordinate zero)?

We show that if such an isomorphism exists, then $A$ and $A'$ have the same cyclic factors; if all cyclic factors have order larger than $2$, the map $f$ is diagonal up to permutation, hence sends coordinate hyperplanes to coordinate hyperplanes.

As a model application, we show using twisted signatures that there exists a family of compact 4-manifolds $X(n)$ with $H_1 X(n) = \mathbb Z/n$ with the property that $\prod X(n_i) \cong \prod X(n'_j)$ if and only if the factors may be identified (up to permutation), and that the induced map on first homology is represented by a diagonal matrix.
\end{abstract}


\end{frontmatter}

\tableofcontents

\section{Introduction}
Let $A = \prod_{i=1}^m \Bbb Z/n_i$ be a product of non-trivial finite cyclic groups. We say that $$Z_i = \{(x_1, \cdots, x_m) \in A \mid x_i = 0\}$$ is the $i$'th coordinate hyperplane (the set of points where the $i$'th coordinate is zero). Their union $\bigcup_{i=1}^m Z_i$ is the set of elements $x \in A$ for which at least one coordinate $x_i$ is zero, which we denote $\mathcal Z(A)$.

Taking $Z_i$ as a model example of a hyperplane, we define a hyperplane in an abelian group to be a subgroup with non-trivial cyclic quotient. We will also have occasion to consider affine hyperplanes, or translates of hyperplanes by some $x \in A$. Our first interest is the following informal question: \emph{When can one recover a set of hyperplanes from their union?}

\begin{figure}[h]
\centering
  \includegraphics[width=0.7\linewidth]{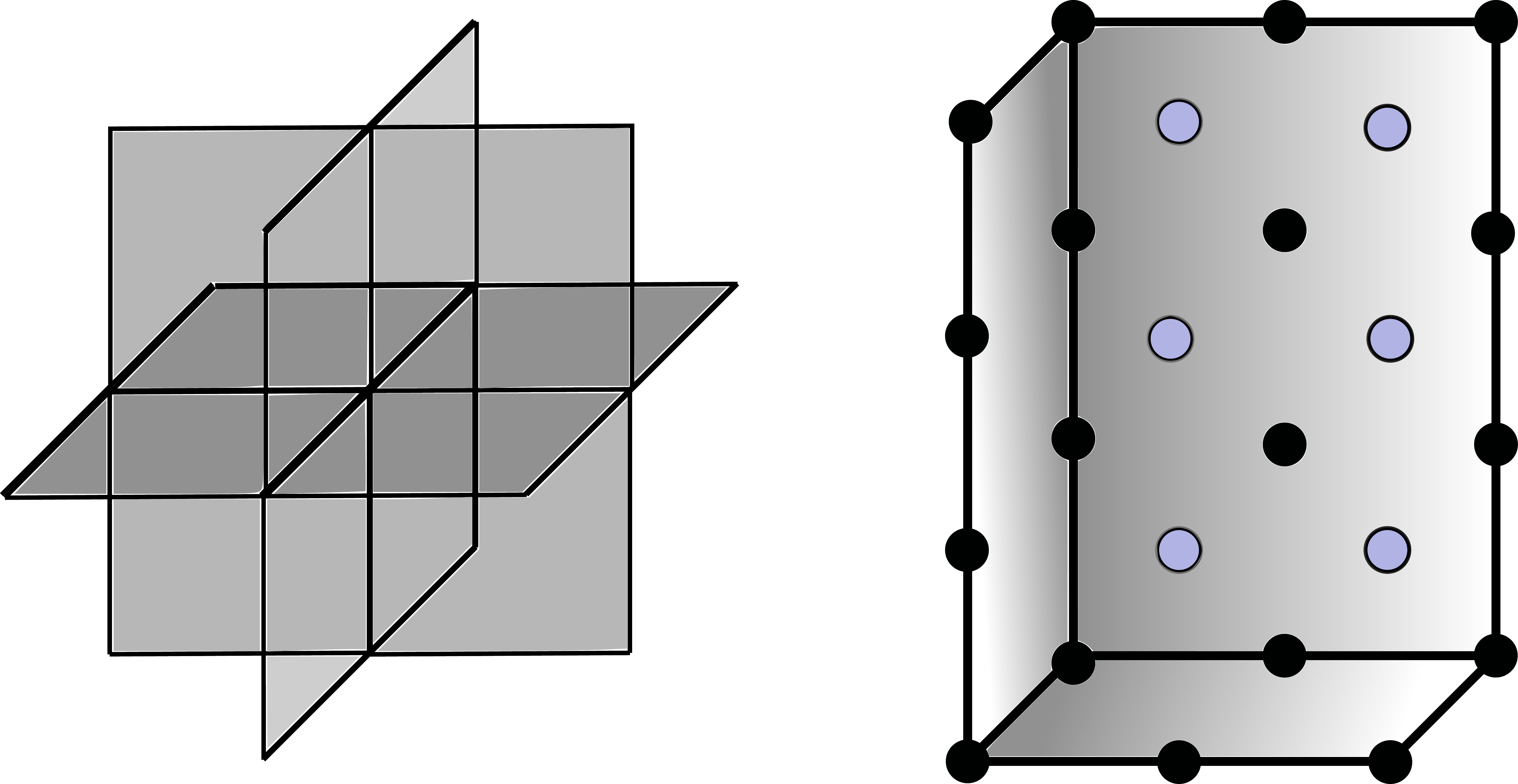}
  \caption{$\mathcal Z(\mathbb R^3)$ and $\mathcal Z(\mathbb Z/2 \times \mathbb Z/3 \times \mathbb Z/4)$. In the latter, we color in the points on $\mathcal Z$ but leave gray the points in its complement.}
  \label{fig:Z-example}
\end{figure}

When $A$ is a vector space, the corresponding question is relatively straightforward. It is visually clear that the only hyperplanes contained in $\mathcal Z(\Bbb R^m)$ are the coordinate hyperplanes. Some linear algebra shows that the same is true for $\mathcal Z\big((\mathbb Z/p)^m\big)$ for $p > 2$.

However, the corresponding claim is false for $A = (\Bbb Z/2)^m$ as soon as $m \ge 3$ (or $m \ge 2$ if we allow affine hyperplanes): we have a containment of the affine line$$\{(x,y) \mid x+y = 1\} = \{(1,0), (0,1)\} \subset \mathcal Z\big ((\mathbb Z/2)^2\big),$$ and this hyperplane is maximal with respect to inclusion among hyperplanes contained in $\mathcal Z(A)$, but is certainly not one of the coordinate hyperplanes. 

\begin{figure}[h]
\centering
  \includegraphics[width=0.65\linewidth]{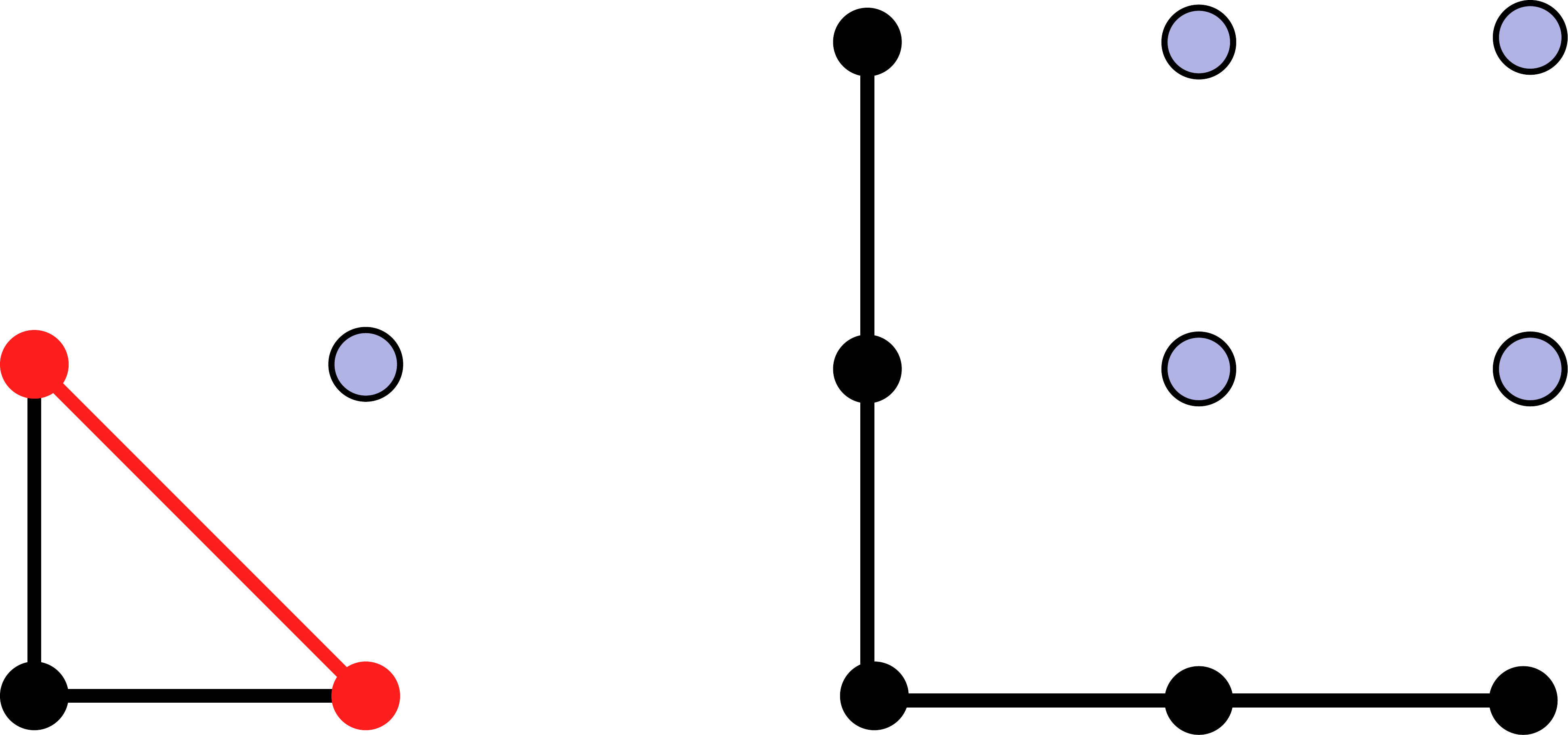}
  \caption{The only affine lines contained in $\mathcal Z\big ((\mathbb Z/3)^2\big)$ are the coordinate lines. By contrast, $\mathcal Z\big ((\mathbb Z/2)^2\big)$ contains the affine line $\{(1,0),(0,1)\}.$}
  \label{fig:counterexample}
\end{figure}

If we restrict ourselves to linear hyperplanes, the first counterexample is the hyperplane
$$\{(x,y,z): x+y+z = 0\} = \{(0,0,0), (1,1,0), (1,0,1), (0,1,1)\} \subset \mathcal Z\big ((\mathbb Z/2)^3\big).$$

The question is further complicated when $A$ is not a vector space, as our definition allows for nested hyperplanes; for instance, $2\Bbb Z/30 \subset 6\Bbb Z/30$ are both non-trivial hyperplanes in $\Bbb Z/30$, as both have cyclic quotient.

To attempt to recover the set of coordinate hyperplanes from their union, we should determine what's special about this set. We propose in Definition \ref{def:hyp-split} the notion of a \emph{hyperplane splitting}, a set $\{H_1, \cdots, H_m\}$ of hyperplanes so that the natural map $\pi: A \to \prod_{i=1}^m A/H_i$ is an isomorphism. Our first result is that --- provided one assumes there are no $\Bbb Z/2$ factors --- one can indeed recover a hyperplane splitting from its union.

\begin{theorem}\label{thm:main1}
Let $A = \prod_{i=1}^m \Bbb Z/n_i$ be a product of finite cyclic groups, with $n_i > 2$ for all $i$. If $\{H_1, \cdots, H_{m'}\}$ is any hyperplane splitting of $A$ with union $\bigcup_{i=1}^{m'} H_i = \mathcal Z(A)$, then $m = m'$ and the $H_i$ are the coordinate hyperplanes of $A$.
\end{theorem}

A more general version of this statement allowing for $n_i = 2$ is proved as Theorem \ref{thm:main1-upgraded} below. One can recover the number of coordinate hyperplanes, and those of odd order quotient are necessarily coordinate hyperplanes. On the other hand, the even order $H_i$ may be more complicated `nearly-coordinate hyperplanes' (which can only appear when at least one factor is $\Bbb Z/2$) modeled on the examples of interesting non-coordinate hyperplanes in $\mathcal Z\big((\Bbb Z/2)^m\big)$ above.\\

As a corollary of our more general statement on hyperplane splittings with given union, we have the following theorem on isomorphisms which preserve the union of coordinate hyperplanes.

\begin{theorem}\label{thm:main2}
Suppose $$A = (\Bbb Z/2)^\ell \times \prod_{i=\ell+1}^m \Bbb Z/n_i \text{ and } A' = (\Bbb Z/2)^{\ell'} \times \prod_{j=\ell'+1}^{m'} \Bbb Z/n'_j,$$ where $n_i, n'_j > 2$ for $i > \ell$ and $j > \ell'$. If $f: A \to A'$ is an isomorphism with $f\big(\mathcal Z(A)\big) = \mathcal Z(A')$, then we have $\ell = \ell', \; m = m'$; additionally, up to a permutation of the $n'_j$, we have $n_i = n'_i$ and a matrix representative of $f$ takes the form $$f = \begin{pmatrix} f_{11} & 0 \\ f_{21} & D \end{pmatrix},$$ where $f_{11}$ is an isomorphism and $D$ is a diagonal matrix.
\end{theorem}

We prove this by examining a special family of hyperplanes in $\mathcal Z(A)$, those containing an element with all but one coordinate equal to $1$. Such hyperplanes are particularly tractable, and classifying them reduces to the case of $p$-groups for $p$ a prime and one factor cyclic of order $p$. When $p$ is odd, such hyperplanes are coordinate hyperplanes; when $p = 2$, the aforementioned nearly-coordinate hyperplanes make an unavoidable appearance.\\

Our initial interest was driven in part by cancellation phenomena for certain types of manifold invariants. An element $x$ of a monoid is \textit{cancellable} if $xy = xz \implies y = z$, and non-cancellable otherwise. 

Suppose one has a collection of mathematical objects together with a product operation $X \boxtimes Y$ (usually Cartesian product or connected sum). We are interested in the question of whether products of these objects are non-cancellative; that is, whether it is possible that $X \boxtimes Y \cong X \boxtimes Y'$ when $Y$ is distinct from $Y'$. 

There is a long history of both cancellation and non-cancellation phenomena for products of manifolds. As an early example of a cancellation phenomenon, the Atiyah--Singer index theorem implies that for 3-dimensional lens spaces $L, L'$ we have $S^1 \times L \cong S^1 \times L'$ if and only if $L \cong L'$ \cite[page 590]{AS-III}. A similar result \cite{KS-1} shows that if $X = T^3 \times P$ and $Y = T^3 \times P'$ for some manifolds $P, P'$ with the homotopy types of spherical space forms, then $T^j \times X \simeq T^j \times Y$ if and only if $X$ is diffeomorphic to $Y$. (See also \cite{BKS,KS-2}.)

Conversely, there are many examples of non-cancellation phenomena, where taking products does lose information. An early example \cite{NC-1} shows that Whitehead's manifold $W$ has $\Bbb R \times \Bbb R^3 \cong \Bbb R \times W$ while $\Bbb R^3 \not\cong W$. Later authors constructed compact manifold examples with few nonvanishing homotopy groups \cite{NC-2}; the introduction of that article points to a great deal of other work constructing non-cancellable products of manifolds.

%For instance, when $X$ and $Y$ are $3$-dimensional lens spaces, $\boxtimes$ is the Cartesian product, and $H$ is the fundamental group, then $X \times X$ and $Y \times Y$ are diffeomorphic if and only if $\pi_1(X) \cong \pi_1(Y)$ \cite[Theorem A]{KS-lens-square}. This ensures non-cancellation when $H$ distinguishes $X$ and $Y$, but allows for cancellation when distinct lens spaces have isomorphic fundamental groups.

We construct infinite families of compact manifolds (with boundary) for which products within this family are cancellative. 

The relevance of Theorem \ref{thm:main2} is as follows. Suppose that each object $X$ is equipped with an abelian group $H(X)$ and a function $\phi_X: H(X) \to \mathbb C$. We say that $\phi$ is \textit{additive} if there is an isomorphism $H(X \boxtimes Y) \cong H(X) \times H(Y)$ so that, with respect to this isomorphism, we have $$\phi_{X \boxtimes Y}(h_1, h_2) = \phi_X(h_1) + \phi_Y(h_2);$$ we say that $\phi$ is \textit{multiplicative} if under this isomorphism we have $$\phi_{X \boxtimes Y}(h_1, h_2) = \phi_X(h_1) \cdot \phi_Y(h_2).$$ The main theorem of \cite{mme:Fourier} shows that when $\phi$ is \textbf{additive}, provided the $\phi_{X_i}$ distinguish the $X_i$ and each has a certain non-vanishing property, one may recover the $\phi_{X_i}$ from the invariant $\phi_{\boxtimes X_i}$. This was applied to classify connected sums of lens spaces up to integer homology cobordism by taking $X_i$ to be 3-dimensional lens spaces, $\boxtimes$ to be connected sum, and $\phi$ to be the Heegaard Floer $d$-invariant. That is, $Y \# L \sim Y \# L'$ implies $L \sim L'$: sums with lens spaces can be cancelled in the homology cobordism monoid.

The argument in the additive case reduces to showing that one can recover the coordinate lines in an abelian group from their union. To make an analogous argument in the multiplicative setting, one is led to precisely Theorem \ref{thm:main2}. To apply our theorem to cancellability of products, we investigate one particular multiplicative invariant: twisted signatures of 4-manifolds.

Associated to a compact oriented 4-manifold $X$ and a homomorphism $\phi: H_1(X) \to S^1$ is its \textit{twisted signature}, an integer $\sigma(X, \phi)$; we review these in Section \ref{sec:application}. The corresponding untwisted signature $\sigma(X, 1)$ is abbreviated to $\sigma(X)$.

We say that a compact connected oriented 4-manifold is \textit{signature-simple} if $\sigma(X, \phi) = 0$ if and only if $\phi$ is the trivial homomorphism. We say a \textit{signature-simple family} is a choice, for each \textbf{odd} integer $n>1$, of a signature-simple manifold $X(n)$ with $H_1(X) \cong \Bbb Z/n$. A (non-constructive) proof that such manifolds exist for each odd $n$ is given in Section \ref{sec:application} below.

\begin{theorem}\label{thm:main3}
Suppose $X(n)$ is a signature-simple family. If there is a homotopy equivalence $$f: \prod_{i=1}^m X(n_i) \to \prod_{j=1}^{m'} X(n_j') \times Y,$$ where $Y$ is a compact connected manifold with $\sigma(Y) \ne 0$ and $H_1(Y) = 0$, then $m = m'$; further, up to a permutation of the $n_j'$, we have $n_i = n_i'$ and $f$ induces a diagonal map on first homology, and $Y$ is a singleton.
\end{theorem}

The set of $\phi \in \text{Hom}(H_1 X, S^1) \cong \prod_{i=1}^m \Bbb Z/m_i$ for which the twisted signatures of the products vanish can be identified with the set of homomorphisms which are trivial on at least one factor. Because a homeomorphism $f$ must preserve these twisted signatures, it follows that it satisfies the assumptions of Theorem \ref{thm:main2}, with all $n_i > 2$.\\

\noindent\textbf{Organization.} In Section 2, we introduce hyperplane splittings and nearly-coordinate hyperplanes, as well as their affine cousins, which appear naturally in our inductive argument. In Section 3, we classify a particular family of affine hyperplanes contained in $\mathcal Z(A)$ (``$\vec z_0$-hyperplanes"), first reducing to the case of $p$-groups --- which are especially amenable to inductive argument --- and then handling the general case. In Section 4, we use our classification of $\vec z_0$-hyperplanes to prove Theorems \ref{thm:main1}-\ref{thm:main2}. Finally, in Section 5 we explain how to apply these results to twisted signatures, and give a construction of signature-simple manifolds.

\section{Affine subgroups and hyperplanes}
In what follows, we will frequently want to refer to terms like $\phi(1, \cdots, 1)$, and sometimes the ambiguous length of the string can cause confusion. We write $\vec 1^m$ for a string of $m$ ones. As is standard, we write $\vec e_i$ for the $i$'th standard basis vector; its length will be clear from context.

In the introduction, we discussed hyperplanes contained in $\mathcal Z(A)$, where $A = \prod_{i=0}^m \Bbb Z/n_i$ is a product of finite abelian groups.\footnote{Note that we now index our products starting at $i = 0$, unlike in the introduction; this choice cleans up the notation in the crucial Definition \ref{def:nearly-coord} of nearly-coordinate hyperplanes.} If $A' = \prod_{i=0}^{m-1} \Bbb Z/n_i$, our inductive approach has us pass from a hyperplane in $A$ to its intersection with $A' \times \{1\}$ projected to $A'$; the new subset $H_1$ then satisfies $H_1 \subset \mathcal Z(A')$. However, because $A' \times \{1\}$ is not a subgroup, neither is $H_1$. Rather, we must investigate the more general notion of \emph{affine hyperplanes}. To start, we should define an affine subgroup. 

\begin{definition}
Let $A$ be an abelian group. An \textbf{affine subgroup} in $A$ is a \textbf{nonempty} subset $S \subset A$ satisfying any of the following three equivalent conditions. 
\begin{enumerate}[label=(\alph*)]
\item $S = x + S'$ is a coset of a subgroup $S' \subset A$. We call $S'$ the \textbf{linear translate} of $S$.
\item $S = f^{-1}(c)$ for some group homomorphism $f: A \to B$.
\item $S$ is closed under affine combinations $\sum_{i=1}^k t_i x_i$, where $x_i \in S$ and $t_i$ are integers with $\sum_{i=1}^k t_i = 1$. 
\end{enumerate}

If $S = S'$ (equivalently, if $0 \in S$), we say that $S$ is a \textbf{linear subgroup}.
\end{definition}

It is straightforward and somewhat pleasant to verify that these conditions are equivalent, and we leave this to the reader. While a `linear subgroup' is by definition a subgroup, the superfluous adjective is useful to contrast with the affine case.

Given an example of an affine subgroup, one can quickly produce more: the image of an affine subgroup under a homomorphism is affine, as is the preimage of an affine subgroup (provided it is nonempty). Similarly, a nonempty intersection of two affine subgroups is again affine.

\begin{definition}
Let $A$ be an abelian group. 
\begin{itemize}
    \item A \textbf{linear hyperplane} in $A$ is a proper subgroup $H \subset A$ with cyclic quotient; equivalently, the zero set of a non-trivial homomorphism $f: A \to \Bbb R/\Bbb Z$. 
    \item An \textbf{affine hyperplane} in $A$ is a coset of a linear hyperplane; equivalently, a nonempty set of the form $f^{-1}(c)$ where $f: A \to \Bbb R/\Bbb Z$ is a non-trivial homomorphism.
\end{itemize}
\end{definition}

We allow for a slightly more general notion of hyperplane splitting than discussed in the introduction, to allow for affine hyperplanes. While this is not relevant in our final application, our arguments apply to the affine context without change.

\begin{definition}\label{def:hyp-split}
A \textbf{hyperplane splitting} of $A$ is an unordered collection $\mathcal H = \{H_0, \cdots, H_m\}$ of affine hyperplanes of $A$ so that the induced map $f: A \to \prod_{i=0}^m A/H'_i$ is an isomorphism. We write $U(\mathcal H) = \bigcup_{i=0}^m H_i$ for the \textbf{hyperplane union}.
\end{definition}

The standard example of a hyperplane splitting is the one introduced at the beginning of this article: $A = \prod_{i=0}^m \Bbb Z/n_i$ has a hyperplane splitting $\{Z_0, \cdots, Z_m\}$ with $Z_j$ the $j$'th coordinate hyperplane, the set of points whose $j$'th coordinate is zero. In fact, given an arbitrary hyperplane splitting of an arbitrary finite abelian group $B$, write $n_i$ for the order of $B/H_i'$; then the homomorphism $$\pi: B \to \prod_{i=0}^m B/H'_i \cong \prod_{i=0}^m \Bbb Z/n_i$$ sends each $H_i$ isomorphically onto a translate of the corresponding coordinate hyperplane. An appropriate translation gives rise to an affine isomorphism $\pi + c: B \cong \prod_{i=0}^m \Bbb Z/n_i$ which sends each $H_i$ to $Z_i$. Thus every hyperplane splitting is isomorphic to the standard splitting by coordinate hyperplanes. 

In a hyperplane splitting $\mathcal H$, the $H_i$ are \textit{maximal} affine hyperplanes contained in $U(\mathcal H)$:

\begin{lemma}\label{lemma:H-maximal}
If $\mathcal H$ is a hyperplane splitting of $A$, the $H_i$ are affine hyperplanes contained in $U(\mathcal H)$, and they are maximal with respect to this property.
\end{lemma}
\begin{proof}
It suffices to verify this for $\mathcal H$ the coordinate hyperplane splitting, where it is clear.
\end{proof}

It would be reasonable to guess that the coordinate hyperplanes are precisely the maximal hyperplanes contained in $\mathcal Z(A)$, or even better, that they are the maximal \textit{affine} hyerplanes contained in $\mathcal Z(A)$ --- so that a maximal affine hyperplane is in fact a linear hyperplane. As discussed in the introduction, some simple computation shows that this is true over $(\Bbb Z/p)^m$ for $p > 2$, but false for $p = 2$, with model example being those vile hyperplanes $$V_m = \{(x_0, \cdots, x_m) \mid x_0 + \cdots + x_m = m\} \subset (\Bbb Z/2)^{m+1}$$ for $m > 0$. The hyperplane $V_m$ lies in $\mathcal Z\big((\Bbb Z/2)^{m+1}\big)$, which is the complement of the single element $\vec 1^{m+1} \in (\Bbb Z/2)^{m+1}$; further, $V_m$ is maximal with respect to this property. 

There are three ways to construct more hyperplanes with properties similar to the $V_m$. \begin{itemize}
    \item One can add a new, independent, factor to form $$V_m \times \Bbb Z/n \subset (\mathbb Z/2)^{m+1} \times \Bbb Z/n.$$
    \item One can \emph{inflate a coordinate} by a factor of $r$ to form the affine hyperplane $$\{(rx_0, \dots, x_m) \mid x_0 + \dots + x_m = m\} \subset (\Bbb Z/2r) \times (\Bbb Z/2)^m$$ with quotient isomorphic to $\Bbb Z/2r$. While this construction can be iterated (performed on different coordinates) and returns another affine hyperplane so long as the inflation factors are pairwise coprime, we will soon restrict attention to a particular family of hyperplanes chosen to exclude any such doubly-inflated examples.
    \item After applying the first construction any number of times and the second construction once, one can permute the coordinates of $A$.
\end{itemize} 

We can rewrite the hyperplane $V_m \times \prod_{i=m+1}^k \mathbb Z/n_i$ with $k$ independent factors as the set of elements satisfying $$x_0 = \phi(x_1,\dots,x_{m+k}) + \phi(\vec 1^{m+k}),$$ where $\phi: (\mathbb Z/2)^m \times \prod_{i=1}^k \mathbb Z/n_i\to \mathbb Z/2$ is given by $\phi(x_1,\dots,x_{m+k}) = x_1 + \dots + x_m$. Taking into account inflation and permutation, this leads us to the following definition.

\begin{definition}\label{def:nearly-coord}
Consider the product of finite cyclic groups $\Bbb Z/2r \times \prod_{i=1}^m \Bbb Z/n_i$, with $r \ge 1$ and $n_i \ge 2$. A \textbf{nearly-$Z_0$}\textbf{ hyperplane} is an affine hyperplane obtained by the following construction. Choose a homomorphism $\phi: \prod_{i=1}^m \Bbb Z/n_i \to \mathbb Z/2$ which is nonzero on $\vec e_i$ only if $n_i = 2$ (we say $\phi$ is supported on $\mathbb Z/2$ factors). Then the corresponding nearly-$Z_0$ hyperplane is
$$Z_\phi = \big\{(x_0,\dots,x_m)\mid x_0 = r\big(\phi(x_1,\dots,x_m) + \phi(\vec 1^m)\big)\big\}$$
We say that $x_0$ is the \textbf{determined coordinate}, and if $\phi(\vec e_i) \neq 0$ we say $x_i$ is a \textbf{determining coordinate}. If $\phi(\vec e_i) = 0$, we say $x_i$ is an \textbf{independent coordinate}.

A hyperplane obtained from $Z_\phi$ by permuting the coordinates such that $x_0$ is sent to $x_j$ is called a \textbf{nearly-$Z_j$}\textbf{ hyperplane}. If a hyperplane is nearly-$Z_j$ for some $j$, we say it is \textbf{nearly-coordinate}.
\end{definition}

\begin{remark}
When $\phi = 0$, we have (as the notation suggests) $Z_\phi = Z_0$. In particular, when all $n_i > 2$, the nearly-coordinate hyperplanes are simply the coordinate hyperplanes.\end{remark}
\begin{remark}
When $r > 1$, the determined coordinate of a nearly-coordinate hyperplane is well-defined: the determined coordinate is the unique coordinate $x_i$ so that the $i$'th projection has $\pi_i(H) \ne \Bbb Z/n_i$. However, when $r = 1$, it is impossible to distinguish the role of the determined coordinate and the determining coordinates; permuting the determined coordinate with some determining coordinate returns the same nearly-coordinate hyperplane.
\end{remark}

The nearly-coordinate hyperplanes can be thought of as hyperplanes where, instead of demanding that $x_i = 0$, we demand that $x_i$ is determined by the values in the other $\Bbb Z/2$-coordinates. The addition of a $\phi(\vec 1^m)$ term is to ensure that $\vec 1^{m+1}$ does not lie in $Z_\phi$, which ultimately guarantees that $Z_\phi \subset \mathcal Z(A)$. In fact, we have better:

\begin{lemma}\label{max:vile}
Any nearly-coordinate hyperplane in $A$ is a hyperplane, is contained in $\mathcal Z(A)$, and is maximal with respect to these properties. 
\end{lemma}
\begin{proof}
The stated properties are invariant under permutations, so it suffices to check them for nearly-$Z_0$ hyperplanes $Z_\phi$. By definition, $Z_\phi$ is the nonempty inverse image of of $r\phi(\vec 1^m)$ under the homomorphism $f_\phi: A \to \Bbb Z/2r$ defined by $$f_\phi(x_0, \cdots, x_m) = x_0 - r\phi(x_1, \cdots, x_m),$$ so $Z_\phi$ is an affine hyperplane. To see that $Z_\phi \subset \mathcal Z(A)$, observe that if all coordinates $x_i$ for $i > 0$ are non-trivial, then $\phi(x_1, \cdots, x_m) = \phi(\vec 1^m)$; the only coordinates which contribute to the value of $\phi$ are the $\Bbb Z/2$ coordinates, where $x_i \ne 0$ means precisely that $x_i = 1$. Thus if $$x_0 = r\left(\phi(x_1,\dots,x_m) + \phi(\vec 1^m)\right)$$ where all $x_i \ne 0$ for $i > 0$, we must have $x_0 = 0$; we have established every element of $Z_\phi$ has some coordinate equal to zero, as desired.

As for maximality, suppose $Z_\phi \subsetneq H$ is properly contained in a larger hyperplane; we will show $H$ contains some element with all entries nonzero. Pick $(y_0, \vec y) \in H \setminus Z_\phi$, meaning that $y_0 \ne r\big(\phi(\vec y + \vec 1^m)\big)$. Observe that there is a distinct vector $(x_0', \vec y) \in Z_\phi \subset H$ with $x_0' = r\big(\phi(\vec y+\vec 1^m)\big)$, and that $(0, \vec 1^m) \in Z_\phi \subset H$ by the definition of $Z_\phi$. Taking affine combinations of these vectors, we see that $$(y_0 - x_0', \vec 1^m) = (y_0, \vec y) - (x_0', \vec y) + (0, \vec 1^m) \in H,$$ which has no component zero, as $y_0 \ne x_0'$.
\end{proof}

Further, at least in a simple case, we can verify that nearly-coordinate hyperplanes are the only affine hyperplanes in $\mathcal Z(A)$.

\begin{lemma}\label{lemma:ord-2}
Every affine hyperplane contained in $\mathcal Z(A)$ with quotient of order $2$ is a nearly-coordinate hyperplane.
\end{lemma}

\begin{proof}
Let $f: A \to \Bbb Z/2$ be a non-trivial homomorphism with $H = f^{-1}(c)$ contained in $\mathcal Z(A)$. Because the claimed properties are invariant under permutation, we suppose that $f$ is supported in precisely the coordinates $0 \le i \le j$, which must have $n_i$ even and $f(\vec e_i) = 1$. In particular, $f(\vec 1^{m+1}) = j+1$. Because $H \subset \mathcal Z(A)$, the element $\vec 1^{m+1} \not\in H$, so $f(\vec 1^{m+1}) \ne c$; we must have $c = j$. 

Now for each $0 \le i \le j$, we have $f(\vec 1^{m+1} + \vec e_i) = j+2 = c$, implying $\vec 1^{m+1} + \vec e_i \in \mathcal Z(A)$. Because this vector has $i$'th coordinate $2$ and all other coordinates $1$, the only coordinate which can be zero is the $i$'th; we must have $n_i = 2$ for all $0 \le i \le j$. Thus \begin{align*}H = f^{-1}(c) &= \{(x_0, \cdots, x_m) \in A\mid x_0 + \cdots + x_j = j\} \\
&= \{(x_0, \cdots, x_m) \in A \mid x_0 = x_1 + \cdots + x_j + j\}.
\end{align*}

This is precisely $Z_\phi$ for $\phi(x_1,\dots,x_m) = x_1 + \cdots + x_j$.
\end{proof}

Before moving on, let us point out that the \textit{hyperplane} hypothesis is significant. 

\begin{example}
The group $(\mathbb Z/n)^n$ contains the affine subgroup whose elements are cyclic permutations of $(0, 1, \cdots, n-1)$, in which all terms have some nonzero coordinate, but for $n>2$ this is not contained in anything comparable to a nearly-coordinate hyperplane. Notice that the quotient of the corresponding linear hyperplane is isomorphic to $(\mathbb Z/n)^{n-1}$, which is cyclic only for $n=2$, where this construction gives the nearly-coordinate hyperplane $V_1$.
\end{example}

The hyperplane condition allows us to rule out these non-examples. Notice that in all of these affine subgroups, \textbf{exactly} one coordinate of each element is zero. We will investigate hyperplanes containing a specific vector with that property: the vector $$\vec z_0 = \vec 1^{m+1} - \vec e_0 = (0, \vec 1^m),$$ which appeared in the proof of Lemma \ref{max:vile}. This will allow us to set up an induction, but it also reduces the number of potential counter-examples in the previous family: 

\begin{example}\label{ex:p-prime}
When $n = p$ is prime, by rescaling all coordinates after the first we may obtain an affine subgroup of $(\Bbb Z/p)^p$ which contains $\vec z_0$ and is contained in $\mathcal Z\big((\Bbb Z/p)^p\big)$. For instance, for $(\Bbb Z/3)^3$, scaling by $(1, 1, 1/2) = (1, 1, 2)$ one obtains the affine subgroup $\{(0, 1, 1), (2, 0, 2), (1, 2, 0)\}$. If we rescale this affine subgroup of $(\Bbb Z/5)^5$ pointwise by the scalars $(1, 1, 1/2, 1/3, 1/4) = (1, 1, 3, 2, 4)$, we obtain the affine subgroup $$\{(0, 1, 1, 1, 1), (4, 0, 3, 4, 2), (3, 4, 0, 1, 3), (2, 3, 2, 0, 4), (1, 2, 4, 3, 0)\}.$$
\end{example}

\section{$\vec z_0$-hyperplanes}
As we will see in Section \ref{sec:main}, our main results will follow from the classification of affine hyperplanes containing the element $\vec z_i = \vec 1^{m+1} - \vec e_i$. Throughout the section, we write $A = \prod_{i=0}^m \Bbb Z/n_i$ unless stated otherwise.

\begin{definition}
A $\vec z_i$\textbf{-subgroup} in $A$ is an affine subgroup $S$ which has $\vec z_i \in S \subset \mathcal Z(A)$. A $\vec z_i$\textbf{-hyperplane} is a $\vec z_i$-subgroup with $A/S'$ cyclic.
\end{definition}

The reader should keep in mind that each $\vec z_i$-subgroup is contained in $\mathcal Z(A)$, even though this property is not mentioned in the name.

If a hyperplane splitting has $U(\mathcal H) = \mathcal Z(A)$, then for each $i$ some hyperplane in the splitting must contain $\vec z_i$. Further, such a vector $\vec z_i$ is unique so long as we avoid order $2$ hyperplanes.

\begin{lemma}\label{lemma:one-h}
Given a $\vec z_i$-subgroup $S$ with $n_i > 2$, no other $\vec z_j$ is contained in $S$.
\end{lemma}
\begin{proof}
Permuting coordinates, it suffices to show that if $S$ is both a $\vec z_0$- and $\vec z_1$-subgroup, then $n_0 = 2$. Supposing $\vec z_0, \vec z_1 \in S$ and taking affine combinations, we see that
$$(2, -1, \vec 1^{m-1}) = 2\vec z_1 - \vec z_0 \in S.$$
Because $S \subset \mathcal Z(A)$, we must have $2 = 0$ in $\Bbb Z/n_0$, and hence $n_0 = 2$. 
\end{proof}

As in the above argument, we may freely permute coordinates and focus our attention on $\vec z_0$-hyperplanes. Our main result in this section, and the main technical result in this paper, is that the \textbf{maximal} $\vec z_0$-hyperplanes are heavily constrained: 

\begin{theorem}\label{thm:general-case}
Every $\vec z_0$-hyperplane is contained in a nearly-$Z_0$ hyperplane. In particular, the \textbf{maximal} $\vec z_0$-hyperplanes are precisely the nearly-$Z_0$ hyperplanes.
\end{theorem}

\noindent We prove this theorem in three steps, each slightly more general than the last. 

First, in Section 3.1, we prove the corresponding claim when $n_0 = p$ is prime and $n_i = p^{j_i}$ for $i > 0$. This is especially tractable, as if one finds that the quotient $A/H'$ contains two $\Bbb Z/n_i$ factors, it is automatically non-cyclic (as these factors are never coprime). In turn, the corresponding version of Theorem \ref{thm:general-case} gives an \textbf{equality} instead of a containment.

Next, in Section 3.2, we handle the case that the first factor is still $\Bbb Z/p$ but the later factors need not be prime powers. We repeatedly use a certain divisibility trick to show that any factor which is not a prime power is necessarily an independent coordinate, thus reducing to the case handled in Section 3.1.

Finally, we prove Theorem \ref{thm:general-case} by analyzing the projection $\pi_0(H)$ to the first coordinate; the previous cases combine to show it must be either trivial (in which case $H$ is contained in $Z_0$) or be equal to $\frac{n_0}{2} \Bbb Z/n_0\Bbb Z$, which allows us to reduce to the already-handled case $n_0 = 2$. 

\subsection{The case of $p$ prime}
We begin with the case where $A$ is a finite abelian $p$-group with first factor $\mathbb Z/p$. %We approach this by inducting on the number of factors of $A$.%: the essential idea is to show that with the addition of a new factor, the only $\vec z_0$-hyperplanes are the nearly-$Z_0$ hyperplanes from the previous factors, with the possible addition of another determining coordinate in the new factor. %This is made simpler by the fact that $p$ is prime, so the projection of a $\vec z_0$-hyperplane to the zeroth coordinate is either $\{0\}$ or $\mathbb Z/p$. 

\begin{proposition}\label{prop:p-group}
Let $A = \Bbb Z/p \times \prod_{i=1}^m \Bbb Z/p^{j_i}$, where $p$ is prime and each $j_i \ge 1$. If $H$ is a $\vec z_0$-hyperplane in $A$, then $H$ is a nearly-$Z_0$ hyperplane $Z_\Phi$.
\end{proposition}

\begin{proof}
We prove this by induction on $m$. The only $\vec z_0$-hyperplane in $\mathcal Z(\mathbb Z/p)$ is $\{0\}$ itself, so the base case $m = 0$ is trivial.

Suppose $H$ is a $\vec z_0$-hyperplane of $A$. Write $A' = \Bbb Z/p \times \prod_{i=1}^{m-1} \Bbb Z/p^{j_i}$, so that $A = A' \times \Bbb Z/p^{j_m}$. Notice that $$H_1 = \{\vec x \in A' \mid (\vec x, 1) \in H\}$$ is a $\vec z_0$-subgroup of $A'$; because there is an embedding $A'/H_1' \hookrightarrow A/H'$ and $A/H'$ is cyclic, so is $A'/H_1'$, and thus $H_1$ is again a $\vec z_0$-hyperplane. By our inductive hypothesis, we have $H_1 = Z_\phi$ for some homomorphism $\phi: \prod_{i=1}^{m-1} \Bbb Z/p^{j_i} \to \Bbb Z/2$. Our goal is to extend $\phi$ to a homomorphism $\Phi: \prod_{i=1}^{m} \Bbb Z/p^{j_i} \to \Bbb Z/2$ and verify that $H = Z_\Phi$. 

As a first step on the way, we claim that the projection to the last factor is $\pi_m(H) = \Bbb Z/p^{j_m}$. Because $\pi_m(H)$ is an affine subgroup of $\Bbb Z/p^{j_m}$ which contains $1$ (as $\pi_m(\vec z_0) = 1$), it takes the form $\pi_m(H) = 1 + p^r \Bbb Z/p^{j_m}$ for some $1 \le p^r \le p^{j_m}$. That $\pi_m(H)$ takes this form means there is some element $(b_0, \vec b, 1 + p^r) \in H$. Let $s_0 = b_0 - \phi(\vec b) - \phi(\vec 1^{m-1})$; by taking affine combinations we see that$$(s_0, \vec 1^{m-1}, 1+p^r) = (b_0, \vec b, 1+p^r) - (\phi(\vec b) + \phi(\vec 1^{m-1}), \vec b, 1) + (0, \vec 1^{m-1}, 1) \in H.$$ 
Because $H \subset \mathcal Z(A)$, at least one coordinate must be zero; we divide into cases based on whether $s_0$ or $1+p^r$ is zero.

\noindent\textbf{Case 1:} Suppose $s_0 = 0$; then $b_0 = \phi(\vec b) + \phi(\vec 1^{m-1})$. By taking affine combinations we see that for all integers $t$ and all $\vec a \in \prod_{i=1}^{m-1} \Bbb Z/n_i$, we have $(\phi(\vec a) + \phi(\vec 1^{m-1}), \vec a, 1 + tp^r)$ equal to the affine combination \begin{align*}t(\phi(\vec b) + \phi(\vec 1^{m-1}), \vec b, 1 + p^r) - t(\phi(\vec b) + \phi(\vec 1^{m-1}), \vec b, 1) + (\phi(\vec a) + \phi(\vec 1^{m-1}), \vec a, 1) \in H.
\end{align*}
This establishes a containment $H_1 \times (1 + p^r\Bbb Z/p^{j_m}) \subset H$. Because $\pi_m(H')$ is equal to $p^r \Bbb Z/p^{j_m}$ and the kernel of $\pi_m$ is given by $H_1'$, comparing cardinalities shows that this containment is an equality. Now $$A/H' \cong A'/H_1' \times \Bbb Z/p^r \cong \Bbb Z/p \times \Bbb Z/p^r;$$ because $A/H'$ is cyclic, we see that $p^r = 1$ and that $H = H_1 \times \Bbb Z/p^{j_m}$. Because $H_1 = Z_\phi$, the hyperplane $H$ is equal to $Z_\Phi$, where $\Phi(\vec x, i) = \phi(\vec x)$.\\

\noindent\textbf {Case 2:} Suppose $s_0 \ne 0$, so that $1+p^r = 0$. For this to be the case, we must have $p = 2$, $r = 0$, and $j_m = 1$. Additionally, $s_0$ must be $1$, so $(1, \vec 1^{m-1}, 0) \in H$. Write $\Phi: \prod_{i=1}^{m-1} \Bbb Z/2^{j_i} \times \Bbb Z/2 \to \Bbb Z/2$ for the function defined by $\Phi(\vec x, i) = \phi(\vec x) + i$; we claim that \begin{align*}H = Z_\Phi &= \{(x_0, \cdots, x_m) \mid x_0 = \Phi(x_1, \cdots, x_m) + \Phi(\vec 1^m)\} \\
&= \{(x_0, \cdots, x_m) \mid x_0 = \phi(x_1, \cdots, x_{m-1}) + \phi(\vec 1^{m-1}) + x_m + 1\}.
\end{align*}
This amounts to the claim that the elements of $H$ take one of the two forms $$(\phi(\vec x) + \phi(\vec 1^{m-1}), \vec x, 1) \quad \quad \text{or} \quad \quad (\phi(\vec x) + \phi(\vec 1^{m-1}) + 1, \vec x, 0).$$
Because $H_1 = Z_\phi$ we know that $(\phi(\vec x) + \phi(\vec 1^{m-1}), \vec x, 1) \in H$; taking affine combinations we see that $$(\phi(\vec x) + \phi(\vec 1^{m-1}) + 1, \vec x, 0) = (\phi(\vec x) + \phi(\vec 1^{m-1}), \vec x, 1) + (1, \vec 1^{m-1}, 0) - \vec z_0  \in H.$$ Thus $Z_\Phi \subset H \subset \mathcal Z(A)$. By Lemma \ref{max:vile}, $Z_\Phi$ is maximal among hyperplanes contained in $\mathcal Z(A)$, so in fact we have $H = Z_\Phi$ as desired.
\end{proof}

\subsection{The general case}
We begin by proving a brief technical lemma which illustrates a divisibility trick used twice in this section.

\begin{lemma}\label{lemma:no-4}
Let $A = \Bbb Z/4 \times A'$, where $A' = \prod_{i=1}^m \Bbb Z/n_i$ is a product of nontrivial cyclic groups. Suppose $S$ is a $\vec z_0$-subgroup. Then the projection to the initial coordinate $\pi_0(S)$ is not equal to $\Bbb Z/4$. 
\end{lemma}
\begin{proof}
Towards a contradiction, suppose $\pi_0(S) = \Bbb Z/4$, so that in particular $1 \in \pi_0(S)$. Thus there exists some $(1, \vec a) \in S$ for some $\vec a \in A$. By taking affine combinations, we see that
$$(k, k(\vec a-\vec 1^m) + \vec 1^m) = k(1, \vec a) + (1-k)(0, \vec 1^m) \in S$$
for each integer $k$. Because $S \subset \mathcal Z(A)$, we know one coordinate must be $0$; in particular, for $k \not\equiv 0 \mod n_0$ we have $k(1-a_i) = 1$ for some $i > 0$, and hence $k$ is invertible mod $n_i$. But if we take $k$ to be the product of all primes dividing $|A|$, then $k \equiv 2 \not\equiv 0 \mod 4$, so every prime dividing $|A|$ (and hence every prime dividing $n_i$) is invertible mod $n_i$. This implies that $n_i = 1$, contradicting our assumption that the cyclic factors are nontrivial.
\end{proof}

We now move to the case where $A$ is a product of cyclic groups with first factor $\mathbb Z/p$. The central idea here is to use the divisibility trick to reduce the problem to the $p$-group case.

\begin{proposition}\label{prop:Z/p}
Let $A = \mathbb Z/p \times A'$, where $A' = \prod_{i=1}^m \mathbb Z/n_i$ is a product of finite cyclic groups. If $H$ is a $\vec z_0$-hyperplane, then $H$ is contained in a nearly-$Z_0$ hyperplane $Z_\Phi$.
\end{proposition}

\begin{proof}
We begin by permuting the coordinates so that $n_i$ is a power of $p$ if and only if $0 \le i \le r$; we write $A = \Bbb Z/p \times A_p \times A_p^\perp$, where $A_p = \prod_{i=1}^r \Bbb Z/n_i$ and $A_p^\perp = \prod_{i=r+1}^m \Bbb Z/n_i$. 

Consider the $\vec z_0$-hyperplane $H_{\vec 1} = \{\vec x \in \Bbb Z/p \times A_p \mid (\vec x, \vec 1^{m-r}) \in H\}$ in $A_p$. By Proposition \ref{prop:p-group}, we have $H_{\vec 1} = Z_\phi$ for some homomorphism $\phi: A_p \to \Bbb Z/2$ supported on $\Bbb Z/2$ coordinates. Writing $\Phi: A_p \times A_p^\perp \to \Bbb Z/2$ for the map $\Phi(\vec x, \vec y) = \phi(\vec x)$, we claim $H \subset Z_\Phi$. 

Towards a contradiction, suppose not. Then there exists some $(y_0, \vec y, \vec w) \in H \setminus Z_\Phi$, which means that $y_0 \ne \phi(\vec y) + \phi(\vec 1^r)$. We know that $(y_0', \vec y, \vec 1^{m-r}) \in Z_\phi \times \vec 1^{m-r} \subset H$, where $y_0' = \phi(\vec y) + \phi(\vec 1^r)$. 

By taking affine combinations we see that $$(k(y_0-y_0'), \vec 1, k(\vec w-\vec 1^{m-r}) + \vec 1^{m-r}) = k(y_0, \vec y, \vec w) -k(y_0', \vec y, \vec 1^{m-r}) + \vec z_0 \in H \subset \mathcal Z(A).$$ Because $y_0 - y_0'$ is nonzero, for any $k \ne 0 \mod p$ the first coordinate is nonzero, so one of the final $m-r$ coordinates must be zero. This means that whenever $k \ne 0 \mod p$, we have $k(1-w_i) \equiv 1 \mod n_i$ for some $r < i \le m$, and in particular $k$ is invertible mod $n_i$. If $k$ is the product of all primes other than $p$ dividing $|A|$, this implies that $n_i$ is a power of $p$; this contradicts the assumption that $n_i$ is only a power of $p$ for $1 \le i \le r$. Thus $H \subset Z_\Phi$, as initially claimed.
\end{proof}

We are now ready to prove Theorem \ref{thm:general-case}. The argument amounts to reducing the first factor to $\mathbb Z/p$ and applying Proposition \ref{prop:Z/p}; the reduction step accounts for the potential for inflation in the initial coordinate of a nearly-$Z_0$ hyperplane.

\begin{proof}[Proof of Theorem \ref{thm:general-case}]
Recall that $A = \prod_{i=0}^m \mathbb Z/n_i$ and $H$ is a $\vec z_0$-hyperplane in $A$. We know that $\pi_0(H)$ is an affine subgroup of $\mathbb Z/n_0$ containing $0$, so it must be a linear subgroup $r\mathbb Z/n_0$ for some $r\mid n_0$. If $r = 0$, then we have $H \subset Z_0$, so we may assume that $\pi_0(H)$ is nontrivial (and hence $0 < r < n_0$).

Let $p\mid n_0/k$ with $p$ prime or $p = 4$. Define
$$\bar H = \{(x_0, \vec x)\mid (n_0x_0/p, \vec x) \in H\} \subset \mathbb Z/p \times \prod_{i=1}^m \Bbb Z/n_i = \bar A.$$
Now $\bar H$ is a $\vec z_0$-hyperplane in $\bar A$ with $\pi_0(\bar H) = \mathbb Z/p$ (in particular, $\bar H$ is not contained in $Z_0$). By Lemma \ref{lemma:no-4} we cannot have $p = 4$, so $p$ is prime; by Proposition \ref{prop:Z/p}, we see that $\bar H \subset Z_\phi$ for a \textbf{nontrivial} homomorphism $\phi: \prod_{i=1}^m \mathbb Z/n_i \to \mathbb Z/2$. In particular, $p = 2$. 

Thus neither $4$ nor an odd prime divides $n_0/r$, so $n_0 = 2r$ and $\pi_0(S) = r\mathbb Z/2r$. This implies that $H$ is contained in the corresponding inflated hyperplane $Z_{\phi} \subset A$, the inflation of $Z_\phi \subset \bar A$ by $r$ in the initial coordinate.
\end{proof}

\section{Proofs of Theorems 1 and 2}\label{sec:main}
Using our results about $\vec z_i$-hyperplanes, we have the following strengthening of Theorem \ref{thm:main1} from the introduction. Because we no longer focus on a particular special coordinate, we return to the convention from the introduction that indexing of products begins at $i=1$.

\begin{theorem}\label{thm:main1-upgraded}
If $\mathcal H = \{H_1, \cdots, H_{m'}\}$ is a hyperplane splitting of $A = \prod_{i=1}^m \Bbb Z/n_i$ with $U(\mathcal H) = \mathcal Z(A)$, then up to a permutation of the $H_i$ we have \begin{enumerate}[label=(\alph*)]
\item When $n_i > 2$, we have $\vec z_i \in H_i$;
\item Each $H_i$ is a nearly-coordinate hyperplane;
\item $m = m'$ and each $A/H_i'$ is a cyclic group of order $n_i$.
\end{enumerate}
\end{theorem}
\begin{proof}
For convenience, we suppose (by permuting the $n_i$) that $n_i > 2$ if and only if $1 \le i \le \ell$. We begin by handling the hyperplanes whose quotients have order larger than $2$. 

Because $U(\mathcal H) = \mathcal Z(A)$ and each $\vec z_i$ lies in $\mathcal Z(A)$, we see that there exists some $H_{j(i)}$ so that $\vec z_i \in H_{j(i)}$.  For $1 \le i \le \ell$, we know from Lemma \ref{lemma:one-h} that if some $H_j$ contains $\vec z_i$, it contains no other $\vec z_{i'}$; for $1 \le i \le \ell$, we see that if $j(i) = j(i')$, we must have $i = i'$. Thus we may permute the $H_i$ so that $j(i) = i$ for all $1 \le i \le \ell$, so that $\vec z_i \in H_i$ for all $i$ with $n_i > 2$. This establishes (a).

Next, observe that when $1 \le i \le \ell$, each $H_i$ is a \textbf{maximal} $\vec z_i$-hyperplane by Lemma \ref{lemma:H-maximal} and the assumption $U(\mathcal H) = \mathcal Z(A)$. By our classification result Theorem \ref{thm:general-case}, we see that when $1 \le i \le \ell$ each $H_i$ is nearly $Z_i$ (and hence has order $n_i$), which establishes (b) and (c) for $1 \le i \le \ell$. 

We conclude by showing that all remaining hyperplanes have order $2$ and that there are $m - \ell$ of them, which immediately gives (c); the rest of (b) then follows from Lemma \ref{lemma:ord-2}.

By definition of hyperplane splitting, we know that $$\pi: A \to \prod_{j=1}^{m'} A/H'_j$$ is an isomorphism, so the corresponding homomorphism $$\bigcap_{i=1}^\ell H'_i \to \prod_{j=\ell+1}^{m'} A/H'_j$$ is also an isomorphism. Now we may identify the domain with $(\Bbb Z/2)^{m-\ell}$ as follows. For $1 \le i \le \ell$ each $H_i$ is nearly $Z_i$, so $H_i'$ takes the form $$H_i' = \{(x_1, \cdots, x_m) \in A \mid x_i = \phi_i(x_{\ell+1}, \cdots, x_m)\},$$ where $\phi_i: (\Bbb Z/2)^{m-\ell} \to \Bbb Z/n_i$ is a homomorphism. Then $$\bigcap_{i=1}^\ell H_i' = \{(x_1, \cdots, x_m) \in A \mid\forall 1 \le i \le \ell \quad x_i = \phi_i(x_{\ell+1}, \cdots, x_m)\}.$$ Thus the projection to the last $m-\ell$ coordinates $$\bigcap_{i=1}^\ell H_i' \to \prod_{j=\ell+1}^m \Bbb Z/2$$ is an isomorphism: the first $\ell$ coordinates are uniquely determined by the last $m-\ell$, which are themselves unconstrained.

We have proved that $\prod_{j=\ell+1}^{m'} A/H'_j$ is isomorphic to $(\Bbb Z/2)^{m-\ell}$. Because each $A/H'_j$ is cyclic and isomorphic to a subgroup of $(\Bbb Z/2)^{m-\ell}$, we see that each $A/H'_j$ is cyclic of order $2$. By comparing cardinalities, we see that $m-\ell = m'-\ell$, as desired.
\end{proof}

Notice that when $n_i > 2$ for all $i$, nearly-coordinate hyperplanes are coordinate, in which case this statement asserts that there is a unique hyperplane splitting (the coordinate hyperplane splitting) with union equal to $\mathcal Z(A)$. This proves Theorem \ref{thm:main1}.

\subsection{Pontryagin duals and Theorem \ref{thm:main2}}
To give the proof of Theorem \ref{thm:main2}, we will want to pass from information about the behavior of hyperplanes in $A$ to information about the behavior of the basis vectors of $A$. The \emph{Pontryagin dual} $A^\vee = \text{Hom}(A, \Bbb R/\Bbb Z)$ is the right tool to do so. For a general finite abelian group $A$, the groups $A$ and $A^\vee$ are noncanonically isomorphic. When we write $A = \prod_{i=1}^m \Bbb Z/n_i$ as a product of finite cyclic groups, there is a canonical isomorphism $$\prod_{i=1}^m \Bbb Z/n_i \cong A^\vee$$ given by sending $\vec x = (x_1, \cdots, x_m)$ to the the homomorphism $\varphi_{\vec x}$, defined by 
$$\varphi_{\vec x}(a_1, \cdots, a_m) = \sum_{i=1}^m \frac{a_i x_i}{n_i}.$$ We freely use this isomorphism in what follows. With respect to this isomorphism, the $j$'th basis vector of $A^\vee$ is given by the homomorphism $$e_j^\vee(a_1, \cdots, a_m) = \frac{a_i}{n_i}.$$
It will be convenient to suppose the factors are ordered so that $$A = (\Bbb Z/2)^\ell \times \prod_{i=\ell+1}^m \Bbb Z/n_i,$$ where $n_i > 2$ for $i > \ell$. We will write $A_2 = (\Bbb Z/2)^\ell \times \{0\}$. Again, $A^\vee$ naturally splits as $A_2^\vee \times \prod_{i=\ell+1}^m \Bbb Z/n_i$.

\begin{lemma}\label{lemma:dual-facts}
The Pontryagin dual group has the following properties.

\begin{enumerate}[label=(\alph*)]
\item An affine hyperplane $H \subset A$ determines a (linear) cyclic subgroup $$C_H = \{\phi \in A^\vee \mid \phi(H') = 0\}$$ with $|C_H| = |A/H'|$.
\item If $f: A \to B$ is a homomorphism, then so is the map $f^\vee: B^\vee \to A^\vee$, defined by $f^\vee(\phi) = \phi f$. If $H \subset B$ is a hyperplane and $J = f^{-1}(H)$, then $f^\vee(C_H) = C_{J}$. 
\item If $A = \prod_{i=1}^m \Bbb Z/n_i$ and $B = \prod_{j=1}^{m'} \Bbb Z/n'_j$, then $f$ can be represented by a matrix $(f_{ij})$. The matrix representation of $f^\vee$ is given by $f^\vee_{ij} = \frac{n_j f_{ji}}{n_i'}$. 
\item If $H$ is a nearly-$Z_i$ hyperplane, then $C_H = \langle e_i^\vee + \phi^\vee\rangle$, where $\phi^\vee \in A_2^\vee$. 
\end{enumerate}
\end{lemma}

We omit the proof of these claims, which are straightforward exercises. These facts in hand, we are ready to complete the proof of Theorem \ref{thm:main2}.

\begin{proof}[Proof of Theorem \ref{thm:main2}]
Suppose $$f: A = (\Bbb Z/2)^\ell \times \prod_{i=\ell+1}^m \Bbb Z/n_i \to (\Bbb Z/2)^{\ell'} \times \prod_{j=\ell'+1}^{m'} \Bbb Z/n'_j = A'$$ is an isomorphism with $f\big(\mathcal Z(A)\big) = \mathcal Z(A')$, and for which $n_i, n'_j > 2$ for all $i > \ell$ and $j > \ell'$.

Write $\mathcal H_A = \{W_1, \cdots, W_m\}$ and $\mathcal H_{A'} = \{Z_1, \cdots, Z_{m'}\}$ for the respective coordinate hyperplane splittings. Then $f$ sends $\mathcal H_A$ to a hyperplane splitting $f(\mathcal H_A)$ with $$U\big(f(\mathcal H_A)\big) = f\big(\mathcal Z(A)\big) = \mathcal Z(A').$$

Write $f(W_i) = \langle \phi_i^\vee\rangle$ for some $\phi_i \in A^\vee$. Because $f$ is an isomorphism, we have $f^{-1}(f(W_i)) = W_i$, so by Lemma \ref{lemma:dual-facts}(b) we see that $f^\vee e_i^\vee = c_i\phi_i^\vee$ for some scalar $c_i$. We need to compute $\phi_i^\vee$.

Applying Theorem \ref{thm:main1-upgraded}, we see that $\ell = \ell'$ and $m=m'$. Each of these are nearly-coordinate; up to a permutation Theorem \ref{thm:main1-upgraded}(a) asserts that for $\ell < i \le m$, each $f(W_i)$ is a $\vec z_i$-hyperplane. By Lemma \ref{lemma:one-h} these combine to show that when $\ell < i \le m$ each $f(W_i)$ is nearly-$Z_i$. Then Lemma \ref{lemma:dual-facts}(d) implies that $\phi_i^\vee = e_i^\vee + \psi_i^\vee$ for some $\psi_i^\vee \in A_2^\vee$ when $\ell < i \le m$. 

This also holds for $1 \le i \le \ell$, though the argument is more roundabout: we know that $f(W_i)$ is nearly-coordinate by Lemma \ref{lemma:ord-2}, so $\phi_i^\vee$ is of the form $c_i(e_j^\vee + \eta_i^\vee)$ for some $e_j^\vee$ (possibly with $i \ne j$) and some $\eta_i^\vee \in A_2^\vee$. Because $e_i^\vee - e_j^\vee$ lies in $A_2^\vee$, so we may simply rewrite this as $c_i(e_i^\vee + \psi_i^\vee)$, taking $\psi_i^\vee = \eta_i^\vee + (e_i^\vee - e_j^\vee)$.

A matrix representation of $f^\vee$ is thus given by $f^\vee = C + M_2$, where $C$ is diagonal with entries $C_{ii} = c_i$ and $M_2$ is supported in the first $\ell$ rows (its $i$'th column is $c_i \psi_j^\vee \in A_2^\vee$). Therefore $f^\vee$ has matrix form $$f^\vee = \begin{pmatrix} f_{11}^\vee & f_{12}^\vee \\
0 & f_{22}^\vee\end{pmatrix},$$ where $f_{22}^\vee$ is diagonal; because $f^\vee$ is an isomorphism and this matrix is block-diagonal, $f_{11}^\vee$ must be an isomorphism. 

By applying Lemma \ref{lemma:dual-facts}(c), we see that the matrix representation of $f$ is a weighted transpose of the matrix representation of $f^\vee$, from which we see that $$f = \begin{pmatrix} f_{11} & 0 \\ f_{21} & D \end{pmatrix},$$ where $f_{11}$ is invertible and $f_{22} = D$ is diagonal.
\end{proof}

\begin{remark}
In fact, because Theorem \ref{thm:main1-upgraded} allows for affine hyperplanes (and affine hyperplane splittings), Theorem \ref{thm:main2} can be slightly strengthened: the map $f(x) = g(x) + c$ can be taken to be an affine isomorphism (so $g$ is an isomorphism and $c \in A'$ is some constant). The conclusion for $g$ is the same, and the conclusion for $c$ is $c \in A_2 \times \{0\}$. 
\end{remark}

\section{An application to topology}\label{sec:application}
Given a compact oriented $2n$-dimensional manifold $X$ (possibly with boundary), consider the cohomology groups $H^n(X;\Bbb C)$ and the relative cohomology groups $H^n(X, \partial X;\Bbb C)$; we suppress the complex coefficients when clear from context. There is a natural map $H^n(X, \partial X) \to H^n(X)$, and intermediate to these lies $$\widehat H^n(X) = \text{image}\left(H^n(X, \partial X) \to H^n(X)\right).$$ The cup product pairing $$H^n(X, \partial X) \otimes H^n(X, \partial X) \to H^{2n}(X, \partial X) \cong \Bbb C$$ descends to bilinear pairings $$H^n(X, \partial X) \otimes H^n(X) \xrightarrow{B_X} \Bbb C \leftarrow H^n(X) \otimes H^n(X, \partial X)$$ which are nondegenerate by Lefschetz duality \cite[Theorem 3.43]{Hatcher}. The existence of these lifts imply that $B_X$ descends to a pairing on $\widehat H^n(X)$, while nondegeneracey of $B_X$ implies $$\widehat B_X: \widehat H^n(X) \otimes \widehat H^n(X) \to \Bbb C$$ is nondegenerate. 

When $X$ is $4k$-dimensional, so that $n = 2k$ and $\widehat B_X$ is a nondegenerate symmetric bilinear pairing, we say that its \textit{signature} is $\sigma(X) = \sigma(\widehat B_X)$, the number of positive eigenvalues of $\widehat B_X$ minus the number of negative eigenvalues (both counted with multiplicity).

When $X$ is not simply connected, one may take (co)homology with \textit{local coefficients} \cite[Section 3.H]{Hatcher}. In particular, given an element $\phi \in \text{Hom}(\pi_1 X, S^1)$, we may construct a Hermitian line bundle $\Bbb C_\phi$ over $X$ with monodromy $\phi$. There is again a nondegenerate Hermitian pairing $H^n(X, \partial X; \Bbb C_\phi) \otimes H^n(X; \Bbb C_\phi) \to \Bbb C$, now given by the composite $$H^n(X, \partial X; \Bbb C_\phi) \otimes H^n(X; \Bbb C_\phi) \to H^{2n}(X, \partial X; \Bbb C_\phi \otimes \Bbb C_\phi) \xrightarrow{(z_1 \otimes z_2) x \mapsto \langle z_1, z_2\rangle x} H^{2n}(X, \partial X; \Bbb C) \cong \Bbb C,$$ where the final map is induced by the Hermitian pairing $\langle -, -\rangle : \Bbb C_\phi \otimes \Bbb C_\phi \to \Bbb C$. 

Once again we consider the intermediate middle cohomology groups $\widehat H^n(X; \Bbb C_\phi)$, and the restricted pairing $\widehat B_{X, \phi}$ is nondegenerate on $\widehat H^n(X; \Bbb C_\phi)$. This is discussed in \cite[Section 2]{APS2}, where one finds the following definition of twisted signatures.

\begin{definition}
If $X$ is a compact oriented $4k$-dimensional manifold, its \textbf{twisted signature} with respect to $\phi \in \textup{Hom}(\pi_1 X, S^1)$ is $$\sigma(X, \phi) = \sigma\left(\widehat B_{X, \phi}\right).$$ 
We say that $X$ is \textbf{signature-simple} if $X$ is connected and $\sigma(X, \phi) = 0$ if and only if $\phi$ is the trivial homomorphism. 
\end{definition}

Signature-simple manifolds are difficult to find (in particular see Remark \ref{rmk:sigsimp} below). It requires, among other things, that each $\widehat H^n(X, \partial X; \Bbb C_\phi)$ is nonzero (for $\phi \ne 1$), which implies that the universal abelian cover of $X$ has $b_n(\hat X) \ge |H_1(X)|-1$, and that $\partial X$ is a nonempty disjoint union of 3-manifolds with nontrivial first homology. 

We will give examples of signature-simple manifolds later. First, we give an algebraic model example to give a sense of what the intersection form of such a manifold might look like.

\begin{example}
Suppose $X$ has $\pi_1(X) = \Bbb Z/n$ and the universal cover of $X$ has $$\widehat H^n(\tilde X, \partial \tilde X) \cong \{\vec z \in \Bbb C^n \mid \sum z_i = 0\}$$ with $\Bbb Z/n$ acting by cyclic permutations of the coordinates, and whose intersection form is the restriction of the standard inner product on $\Bbb C^n$. Write $V$ to denote this $\Bbb Z/n$-representation.

For each $\phi \in \text{Hom}(\Bbb Z/n, S^1)$, the space $\widehat H^n(X, \partial X; \Bbb C_\phi)$ coincides with the $\phi$-eigenspace of $V$. Precisely, writing $\zeta = e^{2\pi i/n}$ and supposing $\phi(1) = \zeta^j$, for $j \ne 0$ this is the subspace $$V_j = \{\vec z \in \Bbb C^n \mid z_{i+1} =  \zeta^j \cdot z_i, \;\; \sum z_i = 0\} = \Bbb C\left\langle \frac{1}{\sqrt n} (1, \zeta^j, \cdots, \zeta^{j(n-1)})\right\rangle,$$ while for $j = 0$ this eigenspace is trivial. 

The intersection form on $V_j \cong \Bbb C$ is then the standard inner product for $j \ne 0$, so that $$\sigma(X, \phi) = \begin{cases} 1 & \phi \ne 1 \\ 0 & \phi = 1 \end{cases}$$ and $X$ is signature-simple.
\end{example}

We have the following basic properties of the twisted signatures.

\begin{lemma}\label{lemma:sigprop}
Twisted signatures of $4k$-dimensional manifolds $(X, \phi)$ have the following properties: 
\begin{enumerate}[label=(\alph*)]
    \item If $f: X \to Y$ is a homotopy equivalence, then $\sigma(X, f^*\psi) = \pm \sigma(Y, \psi)$. 
    \item If $X$ and $Y$ are of the same dimension, then $$\sigma(X \# Y, \phi \# \psi) = \sigma(X, \phi) + \sigma(Y, \psi).$$
    \item If $\phi \in \textup{Hom}(H_1 X, S^1)$ and $\psi \in \textup{Hom}(H_1 Y, S^1)$, where $X$ and $Y$ are compact oriented manifolds of dimension divisible by $4$, we have $$\sigma(X \times Y, \phi \cdot \psi) = \sigma(X, \phi) \cdot \sigma(Y, \psi).$$
    \item If $X$ is a compact oriented 4-manifold, then $\sigma(X, \phi) - \sigma(X) = -\rho(\partial X, \phi|_{\partial X})$, where $\rho(Y,\alpha)$ is the topological invariant of \cite{APS2}.
\end{enumerate}
\end{lemma}
\begin{proof}
The first statement follows immediately from naturality of the cohomology groups and the cup product; $f$ induces an isomorphism between $\widehat H^{2k}(Y; \Bbb C_\psi)$ and $\widehat H^{2k}(X; \Bbb C_{f^* \psi})$ which sends $\widehat B_{Y, \psi}$ to $\pm \widehat B_{X, f^* \phi}$ depending on whether or not $f$ is orientation-preserving, and hence they have the same signature up to sign. 

The second statement follows from the fact that there is a direct sum decomposition $\widehat H^{2k}(X \# Y; \Bbb C_{\phi \# \psi}) \cong \widehat H^{2k}(X; \Bbb C_\phi) \oplus \widehat H^{2k}(Y; \Bbb C_\psi)$ with respect to which $\widehat B_{X \# Y, \phi \# \psi}$ is identified with $\widehat B_{X, \phi} \oplus \widehat B_{Y, \psi}$. 

The third statement is well-known for untwisted signatures: it follows from a combination of the K\"unneth formula $H^{2k+2\ell}(X \times Y;\Bbb C) \cong \bigoplus_{i+j=k+\ell} H^{2i}(X;\Bbb C) \otimes H^{2j}(Y;\Bbb C)$ and a linear algebra lemma \cite[Lemma 1]{sigma-mult}. Because the homomorphism $H_1(X \times Y) \to S^1$ factors as a product of $\phi: H_1(X) \to S^1$ and $\psi: H_1(Y) \to S^1$, an appropriate K\"unneth formula holds for local coefficients, and the linear algebra fact applies without change to show the desired multiplicativity property.

The final statement is a corollary of the Atiyah--Patodi--Singer index theorem, explicitly stated as \cite[Theorem 2.4]{APS2}. 
\end{proof}

\begin{remark}\label{rmk:sigsimp}
The final statement implies that the property of being `signature-simple' is equivalent to the claim that $\sigma(X) = 0$ and $\rho(\partial X; i^* \phi) \ne 0$ for all $\phi \ne 1$, and hence is determined by $\partial X$ and the homomorphism $H_1(\partial X) \to H_1(X)$.

It follows that $H_1(\partial X) \to H_1(X)$ must be surjective, and in particular $\partial X$ is nonempty. Otherwise, the dual map $$i^*: \text{Hom}(\pi_1(X), S^1) = \text{Hom}(H_1(X), S^1) \to \text{Hom}(H_1(\partial X), S^1) = \text{Hom}(\pi_1(\partial X), S^1)$$ is non-injective, and if $\phi$ is a nontrivial element of its kernel, then $$\sigma(X, \phi) - \sigma(X) = -\rho(\partial X; i^* \phi) = -\rho(\partial X; 1) = \sigma(X) - \sigma(X) = 0,$$ so that $\sigma(X) = 0$ implies $\sigma(X, \phi)$ is also zero: $X$ cannot be signature-simple.
\end{remark}

This in hand, we can prove Theorem \ref{thm:main3}. After giving the proof, we give a construction of an infinite family of signature-simple manifolds.

\begin{proof}[Proof of Theorem \ref{thm:main3}]
Recall that $X(n)$ is a \textit{signature-simple family}, a choice of signature-simple manifold $X(n)$ for each odd integer $n \ge 3$, and that $Y$ is a compact connected manifold with $H_1(Y) = 0$ and $\sigma(Y) \ne 0$. Further, we assume we have a homotopy equivalence $$f: X = \prod_{i=1}^m X(n_i) \to \prod_{j=1}^{m'} X(n'_j) \times Y = X'.$$ Our first goal is to show that $m = m'$ and $n_i = n'_i$ up to permutation.

Observe that $H_1(X) \cong \prod_{i=1}^m \Bbb Z/n_i,$ and this induces an isomorphism $$\text{Hom}(\pi_1 X, S^1) \cong \text{Hom}(H_1 X, S^1) \cong  \prod_{i=1}^m \Bbb Z/n_i,$$ as discussed before Lemma \ref{lemma:dual-facts}. With respect to this isomorphism, we have $$\sigma(X, \phi_1, \cdots, \phi_m) = \prod_{i=1}^m \sigma(X(n_i), \phi_i)$$ by Lemma \ref{lemma:sigprop}(c). By the assumption that the $X(n_i)$ are signature-simple, we see that $\sigma(X, \phi_1, \cdots, \phi_m) = 0$ if and only if some $\phi_i$ is trivial. Thus $\sigma(X)^{-1}(0) = \mathcal Z(A')$ is the union of coordinate hyperplanes in $A' = \text{Hom}(H_1 X, S^1) \cong \prod_{i=1}^m \Bbb Z/n_i$. 

Similarly, $H_1(X') \cong \prod_{j=1}^{m'} \Bbb Z/n'_j$, with no factor from $Y$ because $\tilde H_0(Y) = H_1(Y) = 0$. With respect to this decomposition, we have $$\sigma(X', \psi_1, \cdots, \psi_{m'}) = \sigma(Y)\prod_{j=1}^{m'} \sigma(X(n'_j), \psi_j).$$ Because $\sigma(Y) \ne 0$ and the $X(n'_j)$ are signature-simple, we see once again that $\sigma(X')^{-1}(0) = \mathcal Z(A)$ is the union of coordinate hyperplanes in $A = \text{Hom}(H_1(X'), S^1) \cong \prod_{j=1}^{m'} \Bbb Z/n'_j.$

By Lemma \ref{lemma:sigprop}(a), we see that $\sigma(X, f^* \psi) = \pm \sigma(X', \psi)$. In particular, $f^*$ sends $\sigma(X')^{-1}(0) = \mathcal Z(A)$ to $\sigma(X)^{-1}(0) = \mathcal Z(A')$. Because $n_i > 2$ for all $i$, Theorem \ref{thm:main2} applies; we see that $m = m'$ and after a permutation of the $n'_j$ we have $n_i = n'_i$ and $f^*$ diagonal. It also follows that the pushforward map $f_*$ on homology --- being the transpose of $f^*$ with respect to the given direct-sum decomposition --- is also diagonal. 

Finally, we should show that $Y$ is a singleton. Because $m' = m$ and $n'_i = n_i$ up to permutation, we have $X' = X \times Y$. The K\"unneth theorem applied to $X \simeq X \times Y$ implies that $Y$ is acyclic, as both $X$ and $Y$ have finitely-generated homology groups. Because $\sigma(Y) \ne 0$ but $Y$ is acyclic, $Y$ is zero-dimensional; because $Y$ is connected, it is a singleton.
\end{proof}

We conclude by explaining how to construct a family of signature-simple manifolds. The key is Remark \ref{rmk:sigsimp}: if the map $H_1(\partial X) \to H_1(X)$ is surjective, $\sigma(X) = 0$, and $\rho(\partial X)$ is nonvanishing on every nonzero element in the image of $$i^*: \text{Hom}(H_1 X, S^1) \to \text{Hom}(H_1 \partial X, S^1),$$ then $X$ is signature-simple. So we should investigate manifolds whose boundaries are some well-chosen family with sufficiently understood $\rho$ invariants.

The family we choose for this purpose is the lens spaces: 

\begin{definition}
Let $n > 1$ be any integer and $0 < q < n$ be coprime to $n$. We say a compact connected oriented 4-manifold $X$ \textbf{$H_1$-bounds $L(n,q)$} if we have 

\begin{itemize}
    \item For some $k$ there is an oriented diffeomorphism $\partial X \cong \sqcup_{i=1}^k L(n,q)$,
    \item For each boundary component $L(n,q)_j$, the inclusion map $(i_j)_*: H_1(L(n,q)) \to H_1(X)$ is an isomorphism. 
    \item The map $(i_j)_*$ is the same for all $j$. 
\end{itemize}
\end{definition}

Note that the second condition ensure that $H_1(X) \cong \mathbb Z/n$. It then follows that $\sigma(X, \phi) \cong -k\rho(L(n,q), i^*\phi)$, so whether or not $X$ is signature-simple reduces to a question about the $\rho$-invariants of $L(n,q)$. 

\begin{remark}
    Notice that $X$ bounds some number $k$ of copies of $L(n,q)$, not the lens space itself: the latter is impossible. Because the $\rho$ invariants of $L(n,q)$ lie in $\frac{1}{n} \Bbb Z$, if $X$ $H_1$-bounds $L(n,q)$ we must have $n \mid k$, and the smallest we can expect is $k = n$. We achieve this in Lemma \ref{lemma:sig-lens} below.
\end{remark}

Write $\rho_{n,q}(i) = \rho(L(n,q), \phi_i)$ for the function $\Bbb Z/n \to \Bbb Q$ determined by an appropriate isomorphism $\text{Hom}(H_1 L(n,q), S^1) \cong \Bbb Z/n$. The following lemma is valid for any choice of such isomorphism, though the computation we use is in terms of a specific one.

\begin{lemma}\label{lemma:rho-values}
For all $(n,q)$ coprime and $n$ odd, we have $\rho_{n,q}(k) = 0$ if and only if $k \equiv 0 \mod n.$ 
\end{lemma}
\begin{proof}
The relevant computations are over 50 years old, combining the Atiyah--Bott fixed point formula \cite{AB-LF} and formulas relating trigonometric sums to Dedekind sums (see eg \cite[Chapter 6]{HZ}). A recent exposition giving the desired formula is \cite[Section 3.3.4]{thesis:rho}.

Let $n$ be arbitrary and let $q$ be coprime to $n$. If $\rho_{n,q}(k) \equiv 0 \mod 2\Bbb Z$ for some $0 < k < n$, we will show that $n$ is even. By \cite[Equation (3.6)]{thesis:rho}, we have for $0 < k < n$ $$\rho_{n,q}(k) = -\frac{2q}{n} k^2 + 2k -1 + 2 \left\lfloor \frac{kq}{n}\right\rfloor + 4\sum_{j=1}^{k-1} \left\lfloor \frac{jq}{n}\right\rfloor.$$

Taking the expression $\rho_{n,q}(k)$ modulo $2\Bbb Z$ simplifies it to $$\rho_{n,q}(k) \equiv -\frac{2k^2q+n}{n} \mod 2\Bbb Z.$$ If $\rho_{n,q}(k) \equiv 0 \mod 2\Bbb Z$, multiplying through by $n$ we see that $0 \equiv 2k^2 q \equiv -n \mod 2\Bbb Z$; this implies that $n$ must be even.
\end{proof}

\begin{remark}
It is straightforward but slightly more tedious to extend the proof of Lemma \ref{lemma:rho-values} to $n = 4^k m$ and $m$ odd. Lemma \ref{lemma:rho-values} may still be true for $n = 2^{2k+1} m$, but the argument is necessarily more subtle: for instance, we have $\rho_{6,1}(3) = \rho_{8,1}(2) = 2 \equiv 0 \mod 2\Bbb Z$.
\end{remark}

\begin{corollary}\label{cor:sig-lens-implies-sig-simp}
    When $n$ is odd, any compact oriented 4-manifold $X$ which $H_1$-bounds $L(n,q)$ and has $\sigma(X) = 0$ is signature-simple. 
\end{corollary}

Finally, we move on to existence.

\begin{lemma}\label{lemma:sig-lens}
For all coprime $(n,q)$, there exists a 4-manifold $X$ with signature zero which $H_1$-bounds $L(n,q)$; one may even suppose $\partial X \cong \sqcup_{i=1}^n L(n,q)$. 
\end{lemma}

\begin{proof}
We will construct the desired manifold in three steps.

\begin{enumerate}
    \item The first step is the most complicated and most inexplicit: we construct a 4-manifold $X_0$ with boundary $\partial X_0 \cong \sqcup_n L(n,q)$ and so that the homomorphism $\phi: H_1\big(L(n,q)\big) \to \Bbb Z/n$ defined by the composite $H_1(\sqcup_n L(n,q); \Bbb Z) \cong (\Bbb Z/n)^n \xrightarrow{x \mapsto \sum x_i} \Bbb Z/n$ extends to a (necessarily surjective) homomorphism $\Phi_0: H_1(X_0) \to \Bbb Z/n$. As we now explain, this is equivalent to the claim that $n[L(n,q)]$ represents the zero element in a certain group $\Omega_3(B\Bbb Z/n)$.

    For any group $G$ one may define a \emph{classifying space} $BG$ \cite[Section 1.B]{Hatcher} with the property that homotopy classes of maps $X \to BG$ correspond bijectively to homomorphisms $\pi_1(X) \to G$ up to conjugacy. Further, given any space $Y$, one may associate its oriented bordism group $\Omega_d^{SO}(Y)$, whose elements are equivalence classes of pairs $(M, f)$, where $M$ is a closed oriented $d$-manifold and $f: M \to Y$ is a continuous map; we consider $(M,f)$ equivalent to $0$ if there is a compact oriented $(d+1)-$manifold $W$ with $\partial W = M$ and an extension of $f$ to $F: W \to Y$. Addition in this group is given by disjoint union. 

    Combining these, the closed oriented 3-manifold $L(n,q)$ comes equipped with a map $\pi_1 L(n,q) \to \Bbb Z/n$, hence defines a class $[L(n,q)] \in \Omega_3^{SO}(B\Bbb Z/n)$, and $n[L(n,q)]$ is represented by the disjoint union of $n$ copies of $L(n,q)$. We claim above that $\sqcup_n L(n,q)$ bounds a manifold $W$ equipped with a map $\pi_1(X_0) \to \Bbb Z/n$ (and because $\Bbb Z/n$ is abelian, this homomorphism factors through $H_1(X_0)$), and this is equivalent to the claim that $n[L(n,q)] = 0 \in \Omega_3^{SO}(B\Bbb Z/n)$. To conclude, it suffices to observe that $\Omega_3^{SO}(B\Bbb Z/n) \cong \Bbb Z/n$, as follows for instance from Sections 7 and 15 of \cite{ConnorFloyd} or from the more topological arguments of \cite[Lemma 2(1)]{Gordon-G-sig}.

    \item Next we construct a manifold $X_1$ with the same property as in the previous example but so that $\Phi_1: H_1(X_1) \to \Bbb Z/n$ is an isomorphism\footnote{Better yet, the same construction can be used to ensure that $\pi_1(X_1) \cong \Bbb Z/n$.}. To do so, pick a finite set of disjoint embedded loops $\gamma_1, \cdots, \gamma_k$ whose homology classes generate $\ker(\Phi_0)$. Perform surgery on each: delete a tubular neighborhood of each diffeomorphic to $S^1 \times D^3$, and glue in a copy of $D^2 \times S^2$ on each resulting boundary component. A Mayer--Vietoris argument shows that $H_1(X_1)$ can be identified with the quotient of $H_1(X_0)$ by the subgroup generated by the $[\gamma_i]$ --- that is, $\ker(\Phi_0)$ --- and hence $\Phi_1: H_1(X_1) \to \Bbb Z/n$ is an isomorphism. 
    \item To complete the construction, by taking the connected sum of $X_1$ with copies of $\Bbb{CP}^2$ or $\overline{\Bbb{CP}}^2$ as appropriate --- which have signature $\pm 1$ and trivial first homology --- we may obtain a manifold $X$ with the properties above and which also has $\sigma(X) = 0$.
\end{enumerate}

Now the isomorphism $\partial X \cong \sqcup_n L(n,q)$ is chosen so that the map $i: H_1(\partial X) \to H_1(X) \cong \Bbb Z/n$ is the coordinate-summing map on first homology, and hence the map $$\Bbb Z/n \cong \text{Hom}(H_1 X, S^1) \xrightarrow{i^*} \text{Hom}(H_1 \big(\sqcup_n L(n,q)\big), S^1) = \text{Hom}(H_1 L(n,q), S^1)^n \cong (\Bbb Z/n)^n$$ is the diagonal map. %Finally, together with additivity of the $\rho$ invariant under disjoint unions, it follows from Lemma \ref{lemma:sigprop}(d) that $$\sigma(X, \phi) = \sigma(X, \phi) - \sigma(X) \cong -\rho(\sqcup_n L(n,q); \phi|_{\partial X}) = -n \rho(L(n,q), \phi|_{\partial X}).\qedhere$$
\end{proof}

\begin{remark}
    When $n$ is odd, the first step above may be made much more explicit: such a manifold (with $k = n^2$) is described in \cite[Section 5.1]{RRX} inspired by the objects $Q(\theta_1, \theta_2)/\Bbb Z_m$ constructed in \cite[page 167]{Gordon-G-sig}. 
\end{remark}

Corollary \ref{cor:sig-lens-implies-sig-simp} and Lemma \ref{lemma:sig-lens} combine to show

\begin{corollary}
For odd $n$, there exists a signature-simple $4$-manifold $X$ with $H_1(X) \cong \Bbb Z/n$.
\end{corollary}

\section*{Acknowledgements}
Both authors appreciate Siddhi Krishna's time teaching us how to use Inkscape. The first author would like to thank Danny Ruberman for his thoughts on an early version of these ideas (as well as his rightful skepticism about applying them to lens spaces). The second author would like to thank Sophia Fanelle for a careful reading of this paper. Finally, both authors thank the anonymous referee for helpful suggestions. 

\section*{Funding}
While working on this article, the second author was supported by the Columbia University Science Research Fellows Program.

\bibliography{main.bib}
\bibliographystyle{plain}
\end{document}